\documentclass{article}

\pdfoutput=1

\usepackage{a4wide}
\usepackage[all]{xy}
\usepackage{subfigure}
\usepackage{amsmath}
\usepackage{amsfonts}
\usepackage{amssymb}
\usepackage{graphicx}
\usepackage{multirow}

\def\almostsurely{\textit{almost surely}}
\def\almostsure{\textit{almost sure}}

\def\B{\mathrm{B}}

\def\convdist{\stackrel{\textrm{d}}{\rightarrow}}
\def\convprob{\stackrel{\textrm{\footnotesize pr}}{\rightarrow}}

\def\E{\mathrm{E}}

\def\given{\,{\Big |}\,}

\def\iVud{ {\buildrel {\circ \kern 4ex} \over \Vud} }
\def\I{\;\mathbb{I}}

\def\LL{{\cal L}}

\def\N{\mathrm{N}}
\def\OO{\mathrm{O}}

\def\oo{\mathrm{o}}
\def\P{\mathrm{P}}

\newcommand{\quadra}{\sqcup \!\!\!\!\sqcap}

\def\R{\mathbb{R}}

\def\tr{\mathrm{tr}}
\def\u{\mathrm{u}}

\def\Var{\mathrm{Var}}

\newcommand{\Vud}{{\mathrm{V}^{(u,d)}}}

\newtheorem{lemma}{{\bf Lemma}}

\newtheorem{theorem}{{\bf Theorem}}

\newtheorem{remark}{{\bf Comment}}
\newenvironment{demo}{\noindent
\textit{Proof.} }{\nopagebreak\samepage\vspace*{-0.5cm}
\begin{flushright}$\quadra$\end{flushright}}

\renewcommand{\baselinestretch}{1.5}   

\title{
Almost sure convergence and asymptotical normality of a generalization 
of Kesten's stochastic approximation algorithm for
multidimensional case
}

\author{Pedro Cruz \\ pedrocruz@ua.pt \\ Universidade de Aveiro -- Portugal}

\date{14 June, 2005}
\begin{document}

\maketitle


\begin{abstract}
It is shown the almost sure convergence and asymptotical normality of a generalization 
of Kesten's stochastic approximation algorithm for
multidimensional case.

In this generalization, the step  increases or decreases 
if the scalar product of two subsequente increments 
of the estimates is positive or negative.

This rule is intended to accelerate the entrance in the 
`stochastic behaviour' when initial conditions 
cause the algorithm to behave in a `deterministic fashion' for the
starting iterations.
\end{abstract}


\section{Introduction and problem statement}

We consider the problem of finding the stationary 
point $x^*\in\mathbb{R}^n$ of a vector field
$\varphi : \R^n\to \R^n$ using the stochastic approximation
algorithm
\begin{eqnarray}
x_t & = & x_{t-1} - \gamma(s_{t-1}) y_t, \quad t=1,2,\ldots \label{meqal1}\\
s_t & = & {(s_{t-1}+ \u(-y^T_t y_{t-1}))}^+,\quad t=2,3,\ldots \label{meqal2}   
\end{eqnarray}
where
\begin{itemize}
\item $y_t    =  \varphi(x_{t-1}) + \xi_t$, $y_t\in\mathbb{R}^n$ is the $t^{th}$ measure 
of $\varphi$ perturbated by the random vector $\xi_t\in\mathbb{R}^n$;
\item $a^+:=\max\{a,0\}$;
\item $\u$ is a sigmoid function;
\item The random vector $x_0\in\mathbb{R}^n$, and the random variables $s_0$ and $s_1$ 
are initial problem conditions of the algorithm;
\item $x_t\in\mathbb{R}^n$ is the $t^{th}$ approximation to  
the stationary point $x^*\in\mathbb{R}^n$ 
of $\varphi$.
\end{itemize}

We suppose the following assumptions apply.

\textbf{Assumptions B1}
\begin{enumerate}

\item $\{x_0,\xi_1,\xi_2,\ldots,\}$ are mutually independent random vectors 
where vectors $\xi_i$ are identically distributed with mean zero $\E \xi_t=0$
and finite covariance matrix $S_\xi:=\E \, \xi_t \xi_t^T$. We denote ${\cal F}_t$ the $\sigma-$algebra
made by random vectors $\{x_0, \xi_1, \xi_2, \ldots, \xi_t\}$ and random variables
$s_0$ and $s_1$. Assume $s_0$, $s_1$ are mutually independent random variables  from  $\{x_0,\xi_1,\xi_2,\ldots\}$.

\item There exists positive $\Omega$ such that for each open ball 
$I \subset \B(\Omega)$, $\P(\xi_t \in I)>0$.

\item $\E |x_0| < \infty$.

\end{enumerate}

\textbf{Assumptions B2}
\begin{enumerate}

\item $\gamma(s)$ is a monotone decreasing function defined in $[0,+\infty)$ 
so $\gamma(0)$ will denote the maximum value of the step.

\item $\displaystyle \int_0^\infty \gamma(s) ds=\infty$.

\item $\displaystyle \int_0^\infty \gamma^2(s) ds < \infty$.

\end{enumerate}

\textbf{Assumptions B3}
\begin{enumerate}

\item 
There exists a continuous function $V(x) : \mathbb{R}^n \to \mathbb{R}^+$ such that 
\begin{enumerate}

\item $V(x^*) = 0$;

\item $\nabla^2 V(x) \le M$  for each $x$,  $M > 0$ (the largest eigenvalue of $\nabla^2 V(x)$ is less than $M$);

\item $\varphi(x)^T \nabla V(x) > 0$ for each $x \neq x^*$;

\item For each $\gamma^*<\gamma(0)$ and for each $z_0$, the sequence 
\[
  z_t = z_{t-1} - \gamma^* \varphi(z_{t-1})
\]
converges deterministically for the stationary point $x^*$ and verify that  
$\{V(z_t), t=1,2,\ldots\}$ is a monotonous decreasing sequence.

\end{enumerate}

\item There exists positive $R$ and $\beta_0$ such that 
      \[
          \varphi(x)^T \nabla V(x) \ge 
             \frac{1}{2} \gamma(0) \cdot 
             ( \varphi(x)^T M \varphi(x) + \tr(S_\xi M)) + \beta_0 
      \]
      for $|x-x^*| \geq R$. This condition limits the maximum step 
      $\gamma(0)$ and guarantees $\inf_{x \neq x^*} |\varphi(x)| >0$.


%
%

\end{enumerate}

\textbf{Assumptions B4}

\begin{enumerate}

\item $\u$ is a monotone, increasing and bounded function $\mathbb{R}\to\mathbb{R}$, for which  
\[
\u_+ = \lim_{x\to+\infty} \u(x)>0 \textrm{ e } \u_- = \lim_{x\to-\infty} \u(x).
\]

\item Denote $\E_\omega = \E[\u(X^{(\omega)})]$ where
  \[
  X^{(\omega)}=\inf_{  {|\varphi_1| \le \omega}\atop {|\varphi_2| \le \omega} } [-(\xi_1+\varphi_1)^T(\xi_2+\varphi_2)]\,.
  \]
Define $\E_0 := \lim_{\omega \to 0^+} \E_\omega$. Constant $\E_0$ must be positive.

Figure~\ref{fig:exemplo-u} shows possible example for function $\u$ where cases for
known algorithms are included.

\begin{figure}
\renewcommand{\baselinestretch}{1}   
\centering
\mbox{
  \subfigure[Case of Robbins-Monroe algorithm \cite{a010}.]{
	  \unitlength=1cm
	  \begin{picture}(3,2)(-1.5,-1)
	  \put(-1.5,0){\vector(1,0){3}}
	  \put(0,-1){\vector(0,1){2}}
	  \put(-0.9,0.9){$u(x)$}
	  \put(1.3,-0.3){$x$}
	  \thicklines
	  \put(-1.5,0.7){\line(1,0){3}}
	  \end{picture}
  } \quad
  \subfigure[Case of Kesten algorithm \cite{a008}.]{
	  \unitlength=1cm
	  \begin{picture}(3,2)(-1.5,-1)
	  \put(-1.5,0){\vector(1,0){3}}
	  \put(0,-1){\vector(0,1){2}}
	  \put(-0.9,0.9){$u(x)$}
	  \put(1.3,-0.3){$x$}
	  \thicklines
	  \put(-1.5,0){\line(1,0){1.5}}
	  \put(0,0.7){\line(1,0){1.5}}
	  \end{picture}
  }	  
  \subfigure[Case similiar to Plakhov-Almeida algorithm \cite{a012}.]{
	  \unitlength=1cm
	  \begin{picture}(3,2)(-1.5,-1)
	  \put(-1.5,0){\vector(1,0){3}}
	  \put(0,-1){\vector(0,1){2}}
	  \put(-0.9,0.9){$u(x)$}
	  \put(1.3,-0.3){$x$}
	  \thicklines
	  \put(-1.5,-0.5){\line(1,0){1.5}}
	  \put(0,0.7){\line(1,0){1.5}}
	  \end{picture}
  }	  
  \subfigure[Some generic case.\label{figugeral}]{
	  \unitlength=1cm
	  \begin{picture}(3,2)(-1.5,-1)
	  \put(-1.5,0){\vector(1,0){3}}
	  \put(0,-1){\vector(0,1){2}}
	  \put(-0.9,0.9){$u(x)$}
	  \put(1.3,-0.3){$x$}
	  \thicklines
	  \qbezier(-1.5,-0.5)(0,-0.5)(0,0)
	  \qbezier(0,0)(0,0.7)(1.5,0.7)
	  \end{picture}
  }	  
}
\caption{Examples of function $u$.}\label{fig:exemplo-u}
\end{figure}
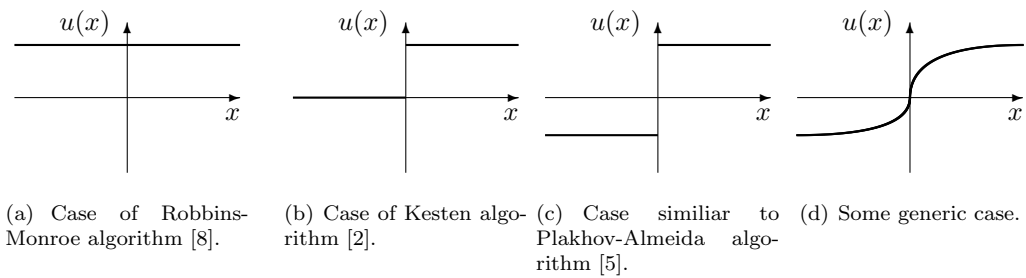

\end{enumerate}

\begin{remark}\label{rem1} 
Suppose we are observing the process (\ref{meqal1}), (\ref{meqal2})
starting in $t_0>1$. This new process, with initial conditions 
$x_{t_0}$, $s_{t_0}$, $s_{t_0+1}$ and the random sequence 
$\xi_{t_0}, \xi_{t_0+1}, \ldots$ also satisfies conditions. Lemma~\ref{c2lemacm4},
for example, makes use of this comment.
\end{remark} 

\begin{remark} 
If $\u$ or the distribution of $\xi_t$
are continuous, then  $\E_0=\E[ \u(-\xi^T_1 \xi_2) ]$.
More, if $\u$ is continuous and verifies
$ \u(x) > -\u(-x)$ when $x \neq 0$,
then B4.2 is valid for any distribution of $\xi_t$ 
with non zero variance.
\end{remark} 

\begin{remark}
We use the following notation for $\varphi$ and $V$: $\varphi'$ denotes a matrix, $\nabla V$ a vector and
$\nabla^2 V$ a matrix.
\end{remark}

\vspace{1cm}

\begin{theorem}\label{c2theorem3}
Suppose Assumptions~B1 to B4 are verified. Then, \almostsurely, 
$\displaystyle \lim_{t\to\infty} x_t = x^*$.
\end{theorem}

\vspace{1cm}

Assumptions for asymptotical normality are all assumptions for \almostsure\ convergence
and three more assumptions: Assumptions~B3.3, B3.4 e B4.3.

\textbf{Assumption B3.3} 
All eigenvalues of $\frac{I}{2} - (1/\E_0) \varphi'(x^*)$ are negative,  
where $I$ is the identity matrix.

\textbf{Assumption B3.4}
Assume Taylor decomposition for $\varphi$,
  \begin{equation}\label{eqa441}
  \frac{|\varphi(x) - \varphi'(x^*) \, (x-x^*)|}{|x-x^*|}=\OO((1),\textrm{ when }x \to x^*\,.
  \end{equation}
\begin{remark}\label{remb34}
From this assumption it follows
  \begin{equation}\label{eqa442}
  \sup |\varphi(x)|/|x-x^*| < \infty
  \end{equation}
because 
  \[
    \frac{|\varphi(x) - \varphi'(x^*) \, (x-x^*)|}{|x-x^*|} \ge
    \frac{|\varphi(x)|}{|x-x^*|}  - |\varphi'(x^*)|
  \]
and so
  \begin{eqnarray*}
  |o(1)| & \ge &  \frac{|\varphi(x)|}{|x-x^*|}  - |\varphi'(x^*)| \\
  \frac{|\varphi(x)|}{|x-x^*|} & \le & |\varphi'(x^*)| - |o(1)| < \infty
  \end{eqnarray*}
\end{remark}

\textbf{Assumption B4.3}
Assume the Taylor decomposition for function $\u$, $\u(x+\Delta x) = \u(x) + \u'(\theta) \Delta x$ 
for $\theta$ between $x$ and $x+\Delta x$.

\begin{theorem}\label{c2theorem4}
Let $x_t$ be defined by (\ref{meqal1}) and (\ref{meqal2}) for which 
\almostsure\ convergence assumptions can be verified. 
Besides, one can also verify Assumptions~B3.3, B3.4 e B4.3.
If $\gamma(s)=1/s$ then
  \begin{equation}
  \sqrt t (x_t - x^*) \convdist \N(0,V)
  \end{equation}
where $\convdist$ denotes convergence in distribution, and $V$ is a positive definite matrix and 
unique solution of the Lyapunov equation (see Theorem~\ref{theorem:lyapunov} in Section~\ref{resultadosusados})
  \begin{equation}\label{m07_1146}
  \left( \frac{I}{2} - (1/\E_0) \varphi'(x^*) \right) 
  (-V) + (-V) 
  \left( \frac{I}{2} - (1/\E_0) \varphi'(x^*) \right)^T
  = (1/\E_0)^2 S_\xi\,.
  \end{equation}
\end{theorem}
\begin{remark}
The explicit solution of equation (\ref{m07_1146}) is
  \[
  (-V) = - \int_{0}^{\infty} e^{W\cdot t} S e^{W^T\cdot t} dt
  \]
where $W =  \frac{I}{2}  - (1/\E_0) \varphi'(x^*)$, $V$ is positive definite. Demonstration
of this result can be find, for example, in  
Theorem 12.3.3 in Lancaster e Tismenetsky \cite{b002}.
\end{remark}

\section{Proof of \almostsure\ convergence}

Demonstration  of the almost sure convergence follows
the work for the unidimensional case by Plakhov e Cruz (2004) \cite{a082}

Without loss of generality we suppose $x^*=0$ so $\varphi(x^*)=0$.

\begin{lemma}\label{c2lemacm1} 
For each $\epsilon>0$ exists $m=m(\epsilon)$ such that, \almostsurely, it occurs
\textit{(i)} exists $t$ such that $|x_t|<\epsilon$,
or 
\textit{(ii)} exists $t$ such that $|x_t|<R$ and $s_t\le m$. 
(Remember that $R$ is defined in B3.2)
\end{lemma} 

\begin{demo}
Choose $\epsilon>0$ and define the stopping time
  \[
  \tau=\tau(\epsilon,m)=\inf\{t:|x_t|<\epsilon \textrm{ or } (|x_t|<R \textrm{ and } s_t \le m) \}.
  \]
Our aim is to prove that for some $m$ we have $\P(\tau=\infty)=0$.

Consider the sequence $\E_t = \E[V(x_t)\I(t<\tau)]$. 

We introduce the simplified notation $V(x_t)=V_t$, 
$\I(t<\tau)=\I_t$, $\nabla V(x_t)= \nabla_t$,
$\gamma(s_t)=\gamma_t$, and using that $\I_{t} \le \I_{t-1}$, we obtain
  \begin{equation}\label{meqn1}
  \E_{t} - \E_{t-1} = \E[V_{t} I_{t} - V_{t-1} \I_{t-1}] \le \E[(V_{t}-V_{t-1}) \I_{t-1}].
  \end{equation}
Using Taylor expansion 
  \[
  V_{t} = V(x_{t-1} -  \gamma_{t-1} y_t) = 
            V_{t-1} -  \gamma_{t-1} y_{t}^T \nabla_{t-1} +
            \frac{1}{2} \gamma^2_{t-1} y_{t}^T \nabla^2 V_{t-1}(x^\prime) y_t,
  \]
where $x^\prime$ is a point between $x_t$ and $x_{t-1}$.
Replacing $y_t$ for $\varphi_{t-1} + \xi_t$ and, in agreement with B3.1, one obtains
  \begin{equation}\label{meqn2}
  V_{t} - V_{t-1} \le
       - \gamma_{t-1} \varphi_{t-1}^T \nabla_{t-1}   - \gamma_{t-1} \xi_t^T \nabla_{t-1} +
       \frac{1}{2} \gamma^2_{t-1} (\varphi_{t-1}^T M \varphi_{t-1} + \xi_t^T M \xi_t).
  \end{equation}
Using (\ref{meqn1}) and (\ref{meqn2}) and observing that each 
values $\gamma_{t-1}$, $\varphi_{t-1}$, $\I_{t-1}$ is determined by $x_{t-1}$ and $s_{t-1}$
and so, mutually independent of $\xi_t$ (Condition~B1.1), 
  \begin{eqnarray*}
  \lefteqn{\E_{t}-\E_{t-1} \le } \\
  & \le & \E[
       - \gamma_{t-1} \varphi_{t-1}^T \nabla_{t-1}   - \gamma_{t-1} \xi_{t}^T \nabla_{t-1} +
       \frac{1}{2} \gamma^2_{t-1} 
       (\varphi_{t-1}^T M \varphi_{t-1} + \xi_t^T M \xi_t) \I_{t-1}] = 
       \nonumber \\
  & =  & \E[- \gamma_{t-1} \varphi_{t-1}^T \nabla_{t-1}] +
         \E[- \gamma_{t-1} \xi_{t}^T \nabla_{t-1}]  + \\
  &    & \E[\frac{1}{2} \gamma^2_{t-1} (\varphi_{t-1}^T M \varphi_{t-1}) \I_{t-1}] + \\
         \E[\frac{1}{2} \gamma^2_{t-1} \I_{t-1} ] \cdot \E[ \xi_t^T M \xi_t ] 
  \end{eqnarray*}
then using
\begin{itemize}
\item $\E[- \gamma_{t-1} \xi_{t}^T \nabla_{t-1}]=0$;
\item $\E[ \xi_t^T M \xi_t ] \le \tr(S_\xi M)$;
\end{itemize}
we have
  \begin{eqnarray}
  \lefteqn{\E_{t}-\E_{t-1} \le } 
    \le   \E[ - \varphi_{t-1}^T \nabla_{t-1} +
               \frac{1}{2} \gamma_{t-1} (\varphi_{t-1}^T M \varphi_{t-1} + \tr(S_\xi M))) 
               \gamma_{t-1} \I_{t-1} 
             ]\,. \label{meqn3}
  \end{eqnarray}
If $\I_{t-1}=1$, then \textit{(i)} $|x_t|\ge R$, or \textit{(ii)} $|x_t|\ge \epsilon$ and $s_t\ge m$.
In case \textit{(i)}, using B3.2, one obtains
  \begin{equation}
   -\varphi_{t-1}^T \nabla_{t-1} + \frac{1}{2} \gamma_{t-1} (\varphi_{t-1}^T M \varphi_{t-1} + \tr(S_\xi M)) \le - \beta_0\,.
  \end{equation}
In case \textit{(ii)} is valid that $\gamma_t < \gamma(m)$ and 
define $\delta_\epsilon := \inf \{ \varphi(x)^T \nabla V(x),\textrm{ for all } |x| \geq \epsilon\}$.
In this context 
  \begin{eqnarray}
  \lefteqn{-\varphi_{t-1}^T \nabla_{t-1} + \frac{1}{2} \gamma_{t-1} (\varphi_{t-1}^T M \varphi_{t-1} + \tr(S_\xi M)) \le }
   \nonumber \\
   & \quad \le & -\delta_\epsilon + \frac{1}{2} \gamma(m) (\varphi_{t-1}^T M \varphi_{t-1} + \tr(S_\xi M)) := -\beta(\epsilon,m)
  \end{eqnarray}
We choose $m$ such that $\beta(\epsilon,m) > 0$ and denote  
$\beta=\inf\{ \beta_0, \beta(\epsilon,m) \}$. So, in both cases, the expression  
between parentesis in right side of  (\ref{meqn3}) is less than 
$- \beta \cdot \gamma_{t-1} \I_{t-1}$ and so
  \[
  \E_{t} - \E_{t-1} \le - \beta \cdot \E[\gamma_{t-1} \I_{t-1}].
  \]

Using that $s_t \le s_0 + t \u_+$ and $\E \I_t = \P(t<\tau)$
one have 
  \[
  \E_{t}-E_{t-1} \le -\beta\,\gamma(s_0 + t \u_+)\,\P(t<\tau);
  \]
by $\P(j<\tau) \ge \P(t<\tau)$ when $j<t$ and,
using induction argument, 
  \[
   E_t \le E_1 - \beta \P(t<\tau) \sum_{j=0}^{t-1} \gamma(s_0 + j \u_+)\,.
  \]
where $\tilde E_0 := E( V(x_0) \I(0<\nu) ) < \infty$ by Assumption~B1.4.

Function $V$ is positive for $x\neq x^*$, so $E_t\ge0$, and from here it follows
  \[
  \P(t<\tau) < \frac{\tilde E_0}{ \beta \sum_{j=0}^{t-1} \gamma(s_0+j\u_+)}.
  \]
When $t\to\infty$ and using $\sum_{j=0}^\infty \gamma(s_0+j\u_+)=\infty$ 
(inferred from Assumption~B2.2), one can conclude that
$P(\tau=\infty) = 0$.
\end{demo}


\begin{lemma}\label{c2lemacm2} 
For each $\epsilon>0$ and $m>0$ exists $\delta$ positive 
such that if $|x_0|<R$ and $s_0\le m$ then
\[
\P( \textrm{exists } t,  |x_t|<\epsilon) \ge \delta\,. 
\]
\end{lemma} 
\begin{demo}
We consider function $V$ defined in Assumptions~B4. 
Let
  \begin{eqnarray*}
  \bar \epsilon = \inf\{ V(x), |x| \ge \epsilon \},\textrm{ and} \\
  \bar R = \sup\{ V(x), |x| \le R \}
  \end{eqnarray*}
then $|x_0| \le R \Rightarrow V(x_0) \le \bar R$ and $V(x) < \bar\epsilon \Rightarrow |x| < \epsilon$.

We will show that  $V(x_t) < \bar \epsilon$ for some $t$.
Denote $V_t := V(x_t)$ and considering the decomposition
  \[
  V_{t} = V_0 \frac{V_1}{V_0} \frac{V_2}{V_1} \cdots \frac{V_{t}}{V_{t-1}}
  \]
First define the deterministic process with constant step $\rho \le \gamma(0)$
  \[
  z_t=z_{t-1} - \rho \varphi(z_{t-1}),\quad t=1,2,\ldots
  \]
and by Assumption~B3.1, exists $V(\cdot)$ such that $\{V(z_t)\}$ converges
monotonically to zero.
Using Taylor expansion
  \begin{eqnarray*}
  V(z_t) 
  & = & V(z_{t-1}  - \rho \varphi(z_{t-1}) ) = \\
  & = & V(z_{t-1}) - \rho \varphi(z_{t-1})^T \nabla V(z_{t-1}) + \\
  &   &        +  \frac{\rho^2}{2} \varphi(z_{t-1})^T \nabla^2 V(z^\prime) \varphi(z_{t-1}) \\
  & = & V(z_{t-1}) - \rho \times  \\
  &   &    (\varphi(z_{t-1})^T \nabla V(z_{t-1}) - 
            \frac{\rho}{2} \varphi(z_{t-1})^T \nabla^2 V(z^\prime) \varphi(z_{t-1}) )
  \end{eqnarray*}
for a certain vector $z^\prime$ between $z_t$ and $z_{t-1}$.
Define
  \[
    U(z,\rho) :=\frac{1}{V(z)} \times \left(
         \varphi(z)^T \nabla V(z) - 
         \frac{\rho}{2} \varphi(z)^T \nabla^2 V(z^\prime) \varphi(z) \right)
  \]
where $z'$ is a point between $z$ and $z-\rho\varphi(z)$
and, since $V(z_t)$ decreases monotonically, then it is necessary that $U(\cdot,\cdot)>0$. 
Define 
  \[
  \bar U := \inf_{ {\epsilon \le |z| \le R} \atop {\rho \le \gamma(0)} } U(z,\rho)
  \]
where  $\bar U$ is a positive constant because $U(\cdot,\cdot)>0$ in
$\epsilon \le |z| \le R$ and 
$\rho \le \gamma(0)$.

Now, we consider Taylor expansion using the original process
  \begin{eqnarray*}
  V(x_t) 
  & =   & V(x_{t-1} - \gamma(s_{t-1}) \varphi(x_{t-1}) - \gamma(s_{t-1}) \xi_t)) \\
  & =   & V(x_{t-1} - \gamma(s_{t-1}) \varphi(x_{t-1}) ) - \\
  &     &  -\gamma(s_{t-1}) 
            \xi_t^T \nabla V(x_{t-1} - \gamma(s_{t-1}) \varphi(x_{t-1})) +
            \frac{\gamma(s_{t-1})}{2} 
            \xi_t^T \nabla V^2(x'')  \xi_t
  \end{eqnarray*}
and defining $\zeta_t:=|\xi_t|$ we have for the last term 
 \begin{eqnarray*}
    -\gamma(s_{t-1}) \xi_t^T \nabla V(x_{t-1} - \gamma(s_{t-1}) \varphi(x_{t-1})) + 
    \frac{\gamma^2(s_{t-1})}{2} \xi_t^T \nabla^2 V(x'') \xi_t
    & \le & \\
    \gamma(0) \zeta_t |\nabla V(x_{t-1} - \gamma(s_{t-1}) \varphi(x_{t-1}))| + 
    \frac{\gamma^2(0)}{2} \zeta^2_t M 
    & \le & \\
    \zeta_t C_\xi
  \end{eqnarray*}
with the following justification
\begin{enumerate}
\item imposing $\zeta_t < 1$;
\item given $\epsilon \le |x| \le R$ then $x_{t-1}$ and $\varphi(x_{t-1})$ 
are vectors from a closed and limited set and $\gamma(s_{t-1}) \le \gamma(0)$,
so $\nabla V( x_{t-1} - \gamma(s_{t-1}) \varphi(x_{t-1}) )$ could be bounded.
\end{enumerate}
From definition of function $U(\cdot,\cdot)$,
  \[
  V(x_{t}) \le V(x_{t-1})(1 - \gamma(s_{t-1}) \cdot U(x_{t-1},\gamma(s_{t-1}))) + 
                            \zeta_t \cdot C_\xi
  \]
and using $1/V(x) \le 1 / \bar \epsilon$, for $\epsilon \le |x| \le R$,
and that $\gamma(s_{t-1}) > \gamma(m + (t-1) \cdot u_+)$,
  \begin{eqnarray*}
  \frac{V_t}{V_{t-1}} 
  &  =  & 1 - \gamma(s_{t-1}) \cdot \bar U + \zeta_t \cdot C_\xi/\bar \epsilon \le \\
  & \le & 1 - \gamma(m+(t-1)u_+) \cdot \bar U + \zeta_t \cdot C_\xi/\bar \epsilon\,.
  \end{eqnarray*}
Denoting  $G_t := 1 - \gamma(m+(t-1)u_+) \cdot \bar U$ we have $G_t < 1$. 
Divergence of the series $\sum_t \gamma(m+t\cdot u_+)$ implies that the productory
$\prod_{i=1}^{t-1} G_i$ goes to zero. 
Using that $G_t \le \sqrt{G_t} < 1$ one can choose $\zeta_t$ such that
 \begin{equation}\label{eq16051104}
  G_t + \zeta_t \cdot C_\xi / \bar \epsilon \le \sqrt{ G_t } < 1
 \end{equation}
and
  \[
  \frac{V_t}{V_{t-1}} \le \sqrt{ G_t }
  \]
whenever that $\epsilon \le |x_{t-1}| \le R$ and $|\xi_t| < \zeta_t < 1$.
We choose $n$ such that $\bar R \prod_{i=1}^{n-1} \sqrt{G_t} < \bar \epsilon$
and suppose we have $|x_0| < R$, $s_0 \le m$ and $|\xi_t| < \zeta_t$ when
$1 \le t \le n-1$. Then, for some $t \in \{1, \ldots, n\}$, $|x_t|<\epsilon$
with probability superior to 
  \[
  \delta := \P( |\xi_1| < \zeta_1, |\xi_2|<\zeta_2, \ldots, |\xi_n|<\zeta_n),
  \]
since from Assumption~B1.2 $P(\xi_t \in I) >0$, for any $I$.
\end{demo}


\vspace{5mm}
From Lemmas~\ref{c2lemacm1} and \ref{c2lemacm2} we have for each $\epsilon>0$ that exists $\delta>0$ 
such that for arbitrary initial conditions $x_0$, $s_0$, $s_1$
\[
\P(\textrm{for some $t$}, |x_t| < \epsilon) > \delta.
\]
Then, we can choose a positive integer number $n=n(x_0,s_0,s_1)$
such that
\[
\P(\textrm{for some } t\le n,|x_t|<\epsilon) > \delta/2\,.
\]
Denote $\bar p=\sup \P(\textrm{for each }t, |x_t|\ge\epsilon)$, being the supremum 
over all initial conditions $x_0$, $s_0$, $s_1$. Fix $x_0$, $s_0$, $s_1$; 
then
\begin{eqnarray}
\lefteqn{\P(\textrm{for each $t$},|x_t|\ge \epsilon) = } \nonumber \\
&=& \P(\textrm{for each $t>n$},|x_t| \ge \epsilon \given \textrm{for each $t\le n$},|x_t| \ge \epsilon)\cdot
    \P( \textrm{for each $t\le n$},|x_t| \ge \epsilon ) \le \nonumber \\ 
&\le& \bar p\,(1-\delta/2)\label{eqntriangulo}. 
\end{eqnarray}
Taking supremum of the L.S. of (\ref{eqntriangulo}) 
over all triple $(x_0, s_0, s_1)$ and denote it by $\bar p$.
Then, we obtain the inequality $\bar p \le \bar p\,(1-\delta/2)$ from which $\bar p=0$.
So, we obtain the following Lemma

\begin{lemma}\label{c2lemacm3} 
For each $\epsilon>0$, \almostsurely\ exists $t$ such that $|x_t|<\epsilon$.
\end{lemma} 

\vspace{5mm}


\begin{lemma}\label{c2lemacm4}
Choose $\epsilon > 0$ and $\eta>0$.
Then, exists $\epsilon_1>0$ and $\delta>0$  such that  
if $|x_0| < \epsilon_1$ then
  \[
  \P( \textrm{for some  $t$, } |x_t| < \epsilon \textrm{ and } s_t \ge \eta) > \delta\,.
  \]
\end{lemma}
\begin{demo}
Starting by $x_t = x_0 - \sum_{i=1}^{t} \gamma_{i-1} y_i$ and using Taylor expansion, 
  \begin{eqnarray*}
  V(x_t) & =   & V(x_0 - \sum_{i=1}^{t} \gamma_{i-1} y_i) \le \\
	 & \le & V(x_0) + |\nabla V(x_0)| 
                          \sum_{i=1}^{t} \gamma_{i-1} |y_i| \cos(y_i,\nabla V(x_0)) 
	         + C_1 |\sum_{i=1}^{t} \gamma_{i-1} y_i|^2\,.
  \end{eqnarray*}

To guarantee the increase in step counter $s_t$ required by this Lemma
we consider two conical symmetrical sections  
where vectors $y_t$ will stay and where we impose
a maximum and a minimum length for $|y_t|$,
$y_I \le |y_t| \le y_{II}$, with  $y_I$, $y_{II}$ to be defined.
We take $x_0$ as a reference point with gradient $\nabla_0:=\nabla V(x_0)$. 
As we will see, we are interested in limiting the internal product
  \[
  y^T \nabla V(x_0) = |y_t| \cdot |\nabla_0| \cdot \cos(y_t,\nabla_0)
  \]
We choose $y_\textrm{odd}$ belongs to the conical section on
the opposite  side of vector $\nabla_0$ and $y_\textrm{even}$ to the conical section. 
%
%
We choose a value $\theta$ for the internal angle of the cone centrered in vector $\nabla_0$ with 
$\theta$ belonging to $(0,\pi/2)$.
In this case $\cos(y_t,\nabla_0)$ is limited by 
  \begin{eqnarray}
  -1 & \le \cos(y_t,\nabla_0) \le & -\cos(\theta), \quad\textrm{$t$ odd},
      \label{ml402_1608}\\
  \cos(\theta) & \le \cos(y_t,\nabla_0) \le &  1, \quad\textrm{$t$ even}\,. 
      \label{ml402_1609}
  \end{eqnarray}
Using (\ref{ml402_1608}) and (\ref{ml402_1609}) we have
  \begin{eqnarray} 
  -y_{II} \le |y_t| \cos(y_1,\nabla_0) \le   -y_I \cos(\theta),  
     \quad\textrm{odd case,} \label{ml403_1024} \\
  y_I \cos(\theta) \le |y_t| \cos(y_2,\nabla_0) \le  y_{II},
  \quad \textrm{even case.}    \label{ml403_1025}
  \end{eqnarray}
It is possible to show $V(x_t) < \bar \epsilon$ if we prove 
  \begin{eqnarray}
  V(x_0)
    & < & \bar \epsilon/3; \label{ml402_1518} \\
   \left| \sum_{i=1}^{t} \gamma_{i-1} |y_i| |\nabla_0| \cos(y_i,\nabla_0) \right|
    & < & \bar \epsilon/3;  \label{ml402_1519} \\
  C_1 |\sum_{i=1}^{t} \gamma_{i-1} y_i|^2 
    & < & \bar \epsilon/3. \label{ml402_1520} 
  \end{eqnarray}

From (\ref{ml402_1518}) we can estimate $\epsilon_1$ by Assumption~B3.3.

From (\ref{ml402_1520}) we conclude
  \begin{equation}\label{ml404_1114}
  C_1 |\sum_{i=1}^{t} \gamma_{i-1} y_i|^2 \le   C_1 y^2_{II} \sum_{i=1}^{\infty} \gamma^2_{i-1} < \bar \epsilon/3
  \end{equation}
and from where we can choose $y_{II}$ (by Assumption~B2.2 the series is convergent).

Because $y_t$ belongs to symmmetrical conical sections,
  \[
  \u(- y^T_t y_{t-1}) \le \u( y^2_I \cos(\pi-\theta) ) = \u( - y^2_I \cos \theta ),
  \quad t=1,2,\ldots,n-1
  \]
therefore
  \begin{equation} \label{ml404_1106}
  s_t \ge (t-2) \u( - y^2_I \cos \theta ),\quad t=3,4,\ldots,n\,.
  \end{equation}
To satisfy $s_t \ge \eta$ required by this Lemma's statement,
we assume $y_I \ge y_{II}/2$, and 
  \begin{equation}\label{ml404_1107}
  n-2 \ge \frac{\eta}{\u( - (y^2_{II}/4) \cos \theta)}
  \end{equation}
obtained from (\ref{ml404_1106}).

Developing the L.S. of (\ref{ml402_1519}) we have by (\ref{ml403_1024}) and (\ref{ml403_1025}),
  \begin{eqnarray}
  \lefteqn{
        -y_{II}               \sum_{ i=1 \atop \textrm{(odd)}  }^{t} \gamma_{i-1} +
        y_I     \cos(\theta)  \sum_{ i=1 \atop \textrm{(even)} }^{t} \gamma_{i-1} \leq} \nonumber \\
   & \leq & \sum_{i=1}^{t} \gamma_i |y_i|\,|\nabla_0|\,\cos(y_i,\nabla_0) \leq \\
   & \leq & 
   -y_I    \cos(\theta)   \sum_{ i=1 \atop \textrm{(odd)}  }^{t} \gamma_{i-1} + 
    y_{II}             \sum_{ i=1 \atop \textrm{(even)} }^{t} \gamma_{i-1}\,. \nonumber
   \end{eqnarray}
Odd sum is bigger than even sum if we start at $i=1$. So
  \begin{equation} \label{ml403_1826}
  \left| \sum_{i=1}^{t} \gamma_{i-1} |y_i|\,|\nabla_0|\,\cos(y_i,\nabla_0) \right| 
  \le 
  y_{II}            \sum_{ i=1 \atop \textrm{(odd)} }^{t} \gamma_{i-1} 
  - y_I  \cos(\theta) \sum_{ i=1 \atop \textrm{(even)} }^t     \gamma_{i-1} 
  \end{equation}
Using (\ref{ml403_1826}),  Condition~(\ref{ml402_1519}) is satisfied if
  \begin{equation}\label{ml404_1150}
  y_{II}            \sum_{ i=1 \atop \textrm{(odd)} }^{t} \gamma_{i-1} 
  - y_I  \cos(\theta) \sum_{ i=1 \atop \textrm{(even)} }^t     \gamma_{i-1} 
  \le  \bar\epsilon /3
  \end{equation}
where we can choose $y_I\ge y_{II}/2$.

For each iteration $t$ the values of $\varphi(x_t):=\varphi_t$,
$y_I$, $y_{II}$, $\theta$  are known.
Let
  \[
  v_t := \frac{ (\varphi_{t-1} + \xi_t)^T \nabla_0}{|y_t| \cdot |\nabla_0|}
  \]
and the conditions that define the admissible region 
for each random vector $\xi_t$ are
  \begin{equation}
  \begin{array}{l}
  y_I \le |\varphi_{t-1} + \xi_t| \le y_{II} \\
  \pi \le \cos^{-1}( v_t ) \le \pi - \theta, \quad\textrm{$t$ odd}\\
  0   \le \cos^{-1}( v_t ) \le \theta,\quad\textrm{$t$ even}.\\
  \end{array}
  \end{equation}
We define $\delta_1$ as the smallest probability of the regions defined in each iteration
$t=1,\ldots,n$ and define $\delta := \delta_1^n$.
Probability $\delta_1$ is positive by Assumption~B1.3. 
\end{demo}


\vspace{5mm}

From Lemmas~\ref{c2lemacm3} and \ref{c2lemacm4} it follows that for each 
$\epsilon>0$ and $\eta>0$ the probability  that
for some $t$, $|x_t|<\epsilon$ and $s_t \ge \eta$ be greater than a positive $\delta$,
will depend only on $\epsilon$ and $\eta$. 
Repeating the argument of Lemma~\ref{c2lemacm3} we have
\begin{lemma}\label{c2lemacm5} 
For each $\epsilon>0$ and $\eta>0$, \almostsurely\ exists $t$ such that $|x_t| < \epsilon$
and $s_t\ge \eta$. 
\end{lemma} 


\vspace{5mm}

We define the stopping time $\tau(\epsilon)=\inf\{t:|x_t|\ge \epsilon\}$.

\begin{lemma}\label{c2lemacm6}
For each $0 < \theta < \E_0$ exists a constant $\epsilon_0>0$ 
and a sequence $\pi_n$ such that 
$\lim_{n\to\infty} \pi_n=0$ and 
\[
\P(s_t > s_0 + t \theta-n\textrm{ for each } t < \tau(\epsilon_0) ) > 1 - \pi_n.
\]
\end{lemma}
\begin{demo}
We will show that 
  \[
  \P( \textrm{exists } t < \tau(\epsilon_0) \textrm{ such that } 
         s_t \le s_0 + t \theta-n ) \le \pi_n \to 0\,.
  \]
From B4.2 it follows that for some $\omega_0$ positive exists $\E_{\omega_0} > \theta$ where
$\E_{\omega_0} = \E[ u( X^{(\omega_0)} ) ]$ and 
  \begin{equation}\label{eqnb244}
  X^{(\omega_0)} = \inf_{ {|\varphi_{1}|\le \omega_0} \atop {|\varphi_{2}|\le \omega_0} } 
                            [-(\xi_1+\varphi_{1})^T (\xi_{2}+\varphi_{2})].
  \end{equation}
We choose $\epsilon_0$ such that 
  \[
  \sup_{|x|<\epsilon_0} |\varphi(x)| \le \omega_0
  \]
and define the sequence $\{\tilde s_t\}$ by
  \begin{equation}\label{eqnb1}
  \tilde s_0=s_0;\quad\tilde s_{t} = \tilde s_{t-1} + \u(X_t^{(\omega_0)})
  \end{equation}
where
  \begin{equation}\label{eqnb2}
  X_t^{(\omega_0)} = \inf_{ {|\varphi_{t-1}|\le \omega_0} \atop {|\varphi_{t-2}|\le \omega_0} } 
                            [-(\xi_t+\varphi_{t-1})^T (\xi_{t-1}+\varphi_{t-2})].
  \end{equation}
Comparing (\ref{eqnb1}) and (\ref{eqnb2}) with (\ref{meqal2}),
for $t < \tau(\epsilon_0)$, we obtain 
  \begin{equation}\label{eqnb3}
  \tilde s_t \le s_t.
  \end{equation}
From (\ref{eqnb1}) it follows that 
  \begin{equation}\label{eqnb4}
  \tilde s_t - s_0 = t \E_{\omega_0} + \I_t^\textrm{even}+\I_t^\textrm{odd}
  \end{equation}
where
  \[
  \I_t^\textrm{even} = 
      \sum_{ {i=1} \atop { (i\,\textrm{even} ) }  }^t 
      [\u(X_t^{(\omega_0)})-\E_{\omega_0}], 
    \quad
  \I_t^\textrm{odd}  = 
      \sum_{ {i=1} \atop { (i\,\textrm{odd}  ) }  }^t 
      [\u(X_t^{(\omega_0)})-\E_{\omega_0}]
  \]
where $\I_t^\textrm{even}$ and $\I_t^\textrm{odd}$ 
are sums of independent and identically distributed random variables
with mean zero and variance linear with $t$.
\begin{remark}
Both variables 
$\I_t^\textrm{even}$ e $\I_t^\textrm{odd}$ are asymptotical normal
however they are dependent from each others. We use the following argument
to estimate the probability of their sum:
$X+Y<a$ implies $X<a/2$ or $Y<a/2$
where $X$ and $Y$ are random variables and $a$ a real constant.
Then,
  \[
  \P(X+Y<a) \le \P(X<a/2) + \P(Y<a/2) \simeq 2 \P(X<a/2).
  \]
\end{remark}
So, using that $\Var \I_t^\textrm{even} = t\cdot V_{\I_1}$, we have
  \begin{equation}\label{eq1447}
  \P( \I_t^\textrm{even} + \I_t^\textrm{odd} < 2a) \lesssim
  2\P( \I_t^\textrm{even} < a) \le  
  2\Phi(\frac{a}{\sqrt t\sqrt V_{\I_1}}).
  \end{equation}
From the event $s_t \le s_0 + t \theta - n$, 
we know that $\tilde s_t \le s_t$ for $t<\tau(\epsilon_0)$.
It follows
  \begin{eqnarray}
  \tilde s_t 
     & \le & s_0 + t \theta - n \Leftrightarrow \nonumber \\ 
   s_0 + t \E_{\omega_0} + \I_t^\textrm{even}+\I_t^\textrm{odd}  
     & \le & s_0 + t \theta - n \Leftrightarrow \nonumber \\ 
  \I_t^\textrm{even}+\I_t^\textrm{odd}  
     & \le & - t (\E_{\omega_0} - \theta) - n\,. \label{eq1128:0923}
  \end{eqnarray}
\begin{remark}
We will use the following argument, where $\{ X_i, i=1,\ldots \}$ is a sequence of random variables,
  \begin{equation}\label{eq0921}
  \P( \textrm{exists } t < \tau \textrm{ such that } X_t < a) 
    \le \sum_{i=1}^\tau \P(X_i < a)
    \le \sum_{i=1}^\infty \P(X_i < a)\,.
  \end{equation}
\end{remark}
By (\ref{eq1447}), (\ref{eq1128:0923}) and (\ref{eq0921}) it follows
  \begin{eqnarray*}
   \P( \textrm{exists } t < \tau(\epsilon_0) \textrm{ such that } 
         s_t \le s_0 + t \theta-n ) 
  & \le & \\
  \P(  \textrm{exists } t < \tau(\epsilon_0) \textrm{ such that } 
         \I_t^\textrm{even}+\I_t^\textrm{odd}  \le 
        - t (\E_{\omega_0} - \theta) - n) 
  & \le & \\
  \sum_{i=1}^\infty \P(
         \I_i^\textrm{even}+\I_i^\textrm{odd}  \le 
        - i (\E_{\omega_0} - \theta) - n) 
  & \lesssim & \\
  2 \sum_{i=1}^\infty \P(
         \frac{\I_i^\textrm{even}}{\sqrt{i V_I}} \le 
        - \sqrt i \frac{\E_{\omega_0} - \theta}{\sqrt V_I} - \frac{n}{\sqrt{ i V_I}}) 
  & \le & \\
  2 \sum_{i=1}^\infty \Phi(- \sqrt i K_1 - \frac{n}{\sqrt i} K_2) 
  & := & \pi_n
  \end{eqnarray*}
for certain constants $K_1>0$ and $K_2>0$.
Last series is convergent and so
$\pi_n\to0$, then
  \begin{equation*}
  \begin{array}{l}
  \lefteqn{\pi_n := 
      \P( \textrm{exists $t$ such that } \I_t^\textrm{even}+\I_t^\textrm{odd} \le }\\
  \quad \quad \le  -n-t (\E_{\omega_0}-\theta) ) \to 0\textrm{ when } n\to\infty.
  \end{array}
  \end{equation*}
\end{demo}

\vspace{5mm}


Now, choose $\theta$ and $\epsilon_0$ as in Lemma~\ref{c2lemacm6}, and 
arbitrarily positive values $\epsilon < \epsilon_0$ and $n$, and
define the stopping time 
  \[
  \nu=\nu(n,\epsilon)=\inf\{ t: |x_t| \ge \epsilon \textrm{ or }
                                s_t \le s_0 - n + t\theta \}
  \] 
and choose
  $\epsilon_1 > 0$
such that 
  \[
  \sup_{|x|<\epsilon_1} V(x) < \frac{1}{2} \inf_{|x|>\epsilon} V(x).
  \]

\begin{lemma}\label{c2lemacm7}  
Let $|x_0|<\epsilon_1$, so 
  \[
  \P(\nu < \infty)\le K \int_{s_0-n-1}^\infty \gamma^2(s)ds + \pi_n,
  \] 
where $K$ is a constant depending  on $\epsilon$.
\end{lemma}
\begin{demo}
Using (\ref{meqn2}) on Lemma~\ref{c2lemacm1}, 
  \[
  V_{t} - V_{t-1}   \le   - \gamma_{t-1} \varphi_{t-1}^T \nabla V_{t-1} 
                        - \gamma_{t-1} \xi_t^T \nabla V_{t-1} + 1/2 \gamma_{t-1}^2 
			   ( \varphi_{t-1}^T M \varphi_{t-1} + \xi_t^T M \xi_t)
  \]
and let $V_{t} - V_0 \le I^\prime_t + I^{\prime\prime}_t$ where
  \[
  \begin{array}{l}
  \displaystyle I^\prime_t = 
    \left | \sum_{i=1}^{t} \gamma_{i-1} \varphi_{i-1}^T \nabla V_{i-1} + 
             \gamma_{i-1} \xi_i^T \nabla V_{i-1} \right| \\
  \displaystyle I^{\prime\prime}_t = 
    1/2 \sum_{i=1}^{t} \gamma^2_{i-1} (\varphi_{i-1}^T M \varphi_{i-1} + \xi_i^T M \xi_i).
  \end{array}
  \]
Let $\delta := (1/2) \inf_{|x|>\epsilon} V(x)$.
For $|x_t| > \epsilon$ then $V_t - V_0 > \delta$, therefore,
  \[
  I^\prime_t + I^{\prime\prime}_t \ge V_t - V_0 > \delta,
  \]
implying $I^\prime_t > \delta/2$ or $I^{\prime\prime}_t > \delta/2$.
We wish to estimate $\P(\nu < \infty)$. Denote
  \[
  \begin{array}{l}
  P^\prime          = \P( I^\prime_\nu \I(\nu < \infty) > \delta/2) \\
  P^{\prime\prime}  = \P( I^{\prime\prime}_\nu \I(\nu < \infty) > \delta/2) \\
  \end{array}
  \]
and using Lemma~\ref{c2lemacm6},
  \begin{equation}\label{m10_1137}
  \P( \nu < \epsilon ) \le \pi_n + P^\prime + P^{\prime\prime}.
  \end{equation}

Using Markov's inequality (for example, \cite[p. 59]{b016}), $\I^2(\cdot) = \I(\cdot)$,
and $\I(i-1 < \nu < \infty) < \I( i-1 < \nu)$, 
\begin{eqnarray*}
  P^\prime 
  & \le & 
    \frac{4}{\delta^2} 
    \E[ {I^\prime_\nu }^2 \I^2(\nu < \infty) ] = \\
  & =   & 
    \frac{4}{\delta^2} 
    \E \left[ 
          \left( 
              \sum_{i=1}^{\nu-1} \gamma_{i-1} (\varphi_{i-1}^T + \xi_i^T) \nabla V_{i-1}) 
          \right)^2 \cdot
          \I(\nu < \infty) 
       \right] \\
  & =   & 
    \frac{4}{\delta^2} 
    \sum_{i,j=1}^\infty \E [ 
       \gamma_{i-1} (\varphi_{i-1}^T+\xi_i^T) \nabla V_{i-1} \I(i-1 < \nu) \times \\
  &     & \quad \times
       \gamma_{j-1} (\varphi_{j-1}^T+\xi_j^T) \nabla V_{j-1} \I(j-1 < \nu) ].
\end{eqnarray*}
Recall that variables $\gamma_{i-1}$, $V_{i-1}$, $\I(i-1 < \nu)$ and $\xi_i$ are mutually
independent. We conclude that terms with $i\neq j$ are zero.
So,
  \begin{equation}\label{m10_1141}
  P^\prime 
    \le \frac{4}{\delta^2} \sum_{i=1}^\infty 
       \E[ \gamma_{i-1}^2 (\varphi_{i-1}^T \nabla V_{i-1})^2 
                          (\xi_i^T \nabla V_{i-1})^2 \I(i-1 < \nu) ] 
    \le K^\prime \E \sum_{i=1}^{\nu-1} \gamma_{i-1}^2
  \end{equation}	       
where $K^\prime$ is a constant that verifies
  \[
  (4/\delta^2) \cdot \sup_{|x|<\epsilon} (\varphi_{i-1}^T \nabla V_{i-1})^2 \cdot 
               \sup_{|x|<\epsilon} \E[ \xi_i^T \nabla V_{i-1}]^2 < K^\prime.
  \]
Using $\P(X>\delta/2)\le \frac{\E|X|}{2/\delta}$,
  \begin{equation}\label{m10_1142}
  P^{\prime\prime} \le \frac{2}{\delta} (1/2) \E[ \sum_{i=1}^{\nu-1} \gamma_{i-1}^2(\varphi_{i-1}^TM\varphi_{i-1} + \xi_i^T M \xi_i)]
                   \le K^{\prime\prime} \sum_{i=1}^{\nu-1} \gamma_{i-1}^2
  \end{equation}
where $K^{\prime\prime}$ verifies
  \[
  (2/\delta) \sup_{|x|<\epsilon} \varphi_{t-1}^T M \varphi_i + \E \xi_i^T M \xi_i < K^{\prime\prime}
  \]
using $\E \xi \xi^T := S_\xi$.

For $t<\nu$, $s_t > s_0 + t \theta-n$, then $\gamma_t < \gamma(s_0-n+t\theta)$,
and
  \begin{equation}\label{m10_1143}
  \E \left[ \sum_{i=1}^{\nu-1} \gamma_i^2 \right] < \sum_{i=1}^\infty \gamma^2(s_0-n+i\theta) \le
  \frac{1}{\theta} \int_{s_0-n-1}^\infty \gamma^2(s) ds.
  \end{equation}
Taking $K=\theta^{-1}(K^\prime + K^{\prime\prime})$, from (\ref{m10_1137}),
(\ref{m10_1141}), (\ref{m10_1142}) and (\ref{m10_1143})
we obtain Lemma~\ref{c2lemacm7}. 
\end{demo}

\vspace{1cm}

Now, choose positive $\epsilon<\epsilon_0$ and choose $n$
and $\eta$ such that  
$1-\pi_n - K \int_{\eta-n-1}^\infty \gamma^2(s) ds$ $=: \delta$
be  positive.
Choose also $\epsilon_1=\epsilon_1(\epsilon)$ as defined above.
In agreement with Lemmas~\ref{c2lemacm5} and \ref{c2lemacm7}, \almostsurely\ exists $t_0$
such that $|x_{t_0}|<\epsilon_1$, $s_{t_0} \ge \eta$, and the probability 
for all $t\ge t_0$, $|x_t|<\epsilon$ exceeds $\delta$. 

We define the sequence of stopping times 
$\tau_1=1$, 
  \[
  \tau_{i+1}=\inf\{ \tau > \tau_i : 
    |x_\tau| \ge \epsilon,
    \textrm{ and for some } \tau_i \le t < \tau, \, |x_t|<\epsilon_1 
    \textrm{ and } s_t > \eta \},\quad
  i=1,2,\ldots\,.
  \]
We have
\[
\P(\tau_{i+1}=\infty \,|\, \tau_i<\infty) \ge \delta,
\]
from
\[
\P(\tau_{i+1} < \infty) = \P(\tau_{i+1}< \infty \,|\, \tau_i < \infty) \, \P(\tau_i< \infty)
\le (1-\delta) \, \P(\tau_i < \infty).
\]
So, $\P(\tau_i<\infty) \to 0$ quando $i\to \infty$; 
implying that \almostsurely\ $i_0=\sup\{ i: \tau_i < \infty \}$
is finite.

In accordance to Lemma~\ref{c2lemacm5}, \almostsurely\ exists $t_0 \ge \tau_{i_0}$ such that 
$|x_{t_0}| < \epsilon_1$ and $s_{t_0} > \eta$; from here 
we conclude that $|x_t| < \epsilon$ when $t > t_0$. 
Theorem~1 is proved.\hfill$\Box$

\clearpage

\section{Proof of the asymptotical normality}

The central idea of the proof  follows 
the work of Delyon and Juditsky (1993) \cite{a003}.


\begin{lemma}[Delyon e Juditsky \cite{a003}]\label{delyon_at}
Let  $(\nu_t)$ be a random sequence of real numbers such that $\nu_t \to 0$ \almostsurely\
when $t\to\infty$. Then exists a deterministic sequence $(a_t)$ such that 
  \begin{equation}
  a_t \to 0 \quad \textrm{  and  } \quad \nu_t/a_t \to 0\quad \textrm{\almostsurely}.
  \end{equation}  
\end{lemma}

In what follows $o$ and $O$ have the standard deterministic meaning however 
many times they represent stochastic random variables belonging to ${\cal F}_t$
$\sigma-$algebra of events.

\begin{lemma}\label{lemaxixi}
Let $\{z_i,i=1,\ldots\}$ be a sequence of non-negative random variables verifying $z_i\to0$ \almostsurely,
and let $\{|\xi_i|\}$, be a sequence of iid random variables with finite variances. 
Possibly, variables $z_i$ and $\xi_i$ are dependent.
Then
  \[
  \sum_{i=1}^t z_i \, |\xi_i| = o(t)
  \]
\almostsurely.
\end{lemma}
\begin{demo}
From Lemma~\ref{delyon_at} there exists a deterministic sequence $\{a_i\}$ such that $z_i/a_i\to0$ \almostsurely.
Then $0\le z_i(\omega)/a_i < M(\omega)$ for each elementary event $\omega$.
Denote $\zeta_i := |\xi_i|-\mu$ where $\mu := \E(|\xi|)$, so $\E \zeta_i = 0$ and $\Var \zeta_i < \infty$. 

Let $S_t = \sum_{i=1}^t a_i \zeta_i$. Then $S_t/t\to 0$ in probability by Chebychev inequality.  Then, by Levy's Theorem 
(for example, \cite{b029} p.~211) $S_t/t \to 0$ \almostsurely\ because $\{a_i \zeta_i\}$ is a sequence of
independent random variables. 
(The same result using Kronecker Lemma \cite{b029} because $\sum \Var(a_i \zeta_i/i) < \infty$.)

Then $S_t = o(t)$ \almostsurely\ and
  \begin{eqnarray*}
  \left | \sum_{i=1}^t \frac{z_i}{a_i} \cdot a_i \cdot |\xi_i| \right | & \le & 
            M(\omega) \cdot \sum_{i=1}^t a_i \cdot |\xi_i|  \\
   =   M(\omega) \cdot \sum_{i=1}^t (a_i \cdot \zeta_i + a_i \cdot \mu_{|\xi|}) & = &  
            M(\omega) \cdot o(t)  = o(t)\,\textit{almost surely}.
  \end{eqnarray*}
\end{demo}

Recall definition of $\E_0$ in Assumption~B4.2.
\begin{lemma}\label{lema_st_draft}
Let $s_0$ and $s_1$ be random variables which are initial conditions
of the process $\{s_t\}$,
defined in (\ref{meqal2}).
Then
  \begin{equation} \label{03_1225}
  \gamma(s_t) = 1/s_t = \frac{1}{\E_0 t} (1+o_t), \textrm{ \almostsurely}
  \end{equation}
where $o_t$ is a random variable defined in ${\cal F}_t$ and 
for which $\lim_{t \to \infty} o_t = 0$ \almostsurely.
\end{lemma}
\begin{demo}
Assumption~B4.3 permits the decomposition
  \begin{equation}\label{utaylor}
  \begin{array}{l}
  \u(-y_{i-1} y_i) = \u(-(\varphi_{i-2}+\xi_{i-1})^T (\varphi_{i-1} + \xi_i)) = \\
  \quad = \u(-(\varphi_{i-2}+\xi_{i-1})^T (\varphi_{i-1} + \xi_i)) = \\
  \quad = \u(-\varphi_{i-2}^T \varphi_{i-1} -\varphi_{i-2}^T \xi_i - \varphi_{i-1}^T \xi_{i-1} - \xi_{i-1}^T \xi_i) = \\
  \quad = \u( - \xi_{i-1}^T \xi_i ) + 
          \u'(\theta_i) \times 
           \left(
                  -\varphi_{i-2}^T \varphi_{i-1} -\varphi_{i-2}^T \xi_i - \varphi_{i-1}^T \xi_{i-1}
           \right)
  \end{array}
  \end{equation}
where $\theta_i$ is a point between $-y_{i-1}^T y_i$ and $-\xi_{i-1}^T \xi_i$.
We also have that function $\u'$ is limited and $\varphi(x_i) \to 0$ from where, 
by Lemma~\ref{lemaxixi},
  \begin{eqnarray}
  \sum_{i=1}^t \u'(\theta_i)\varphi^T_{i-2} \varphi_{i-1} &=& o(t) \label{1021} \\
  \sum_{i=1}^t \u'(\theta_i)\varphi^T_{i-2} \xi_i & = & o(t) \label{1022} \\
  \sum_{i=1}^t \u'(\theta_i)\varphi^T_{i-1} \xi_{i-1} & = & o(t)\,. \label{1023} 
  \end{eqnarray}
So, we have
\begin{eqnarray*}
  s_t & = & s_0 + s_1 + \sum_{i=1}^t (\u(-y_{i-1}^T y_i) - \u(-\xi_{i-1}^T \xi_i) ) +  \nonumber\\
      &   & + \sum^t_{\textrm{even}} \u(-\xi_{i-1}^T \xi_i) + 
              \sum^t_{\textrm{odd}} \u(-\xi_{i-1}^T \xi_i) \nonumber \\
      & = &   s_0 + s_1 + \Delta U_t + P_t + I_t. \label{m01_2317}
\end{eqnarray*}
By (\ref{1021}), (\ref{1022}) and (\ref{1023})
 \[
 \Delta U_t = \sum_{i=1}^t (\u(-y_{i-1} y_i) - \u(-\xi_{i-1} \xi_i) ) = o(t)\,\textrm{\almostsurely}.
 \]

Each of the sums $P_t$ and $I_t$ is composed of independent terms
of mean $\E_0$ and finite variance.
By the law of iterated logarithm 
  \[
  P_t + I_t = \E_0 t + \oo(\sqrt{t \log \log t})\,.
  \]
Using $\lim_{t\to\infty} s_0/t=0$ \almostsurely, also for $s_1$, we have
  \[
  s_t  = s_0 + s_1 + \E_0 t + t o_t + \oo(\sqrt{t \log \log t}) = (\E_0 + o_t) t,
  \]
\almostsurely.
Then
  \begin{eqnarray*} 
  s_t  & = & ( \E_0 + o_t ) t = 
             \E_0 t \left(\frac{1}{1-\frac{o_t}{\E_0 + o_t}} \right) = \\
       & = & \E_0 t \left(\frac{1}{1+o_t} \right).
  \end{eqnarray*}
\end{demo}


\textbf{Demonstration of Theorem 2}
We choose $x^*=0$. 
From last Section, we have shown the \almostsurely\  convergence of $x_t \to 0$ and 
in Lemma~\ref{lema_st_draft} we shown the mean beahaviour of $s_t = \E_0 t (\frac{1}{1+o_t} )$
where $o_t\to0$ \almostsurely.

By Lemma~\ref{delyon_at} we can conclude that there exists a sequence $(a_t)$ of positive non random numbers 
such that
  \begin{equation}\label{01_1800}
  a_t \to 0 \quad \textrm{ and } \quad |o_t|/a_t \to 0, \quad |x_t|/a_t \to 0 \quad \textrm{\almostsurely}.
  \end{equation}
\begin{remark}
We provide an explanantion for the above fact. We can make
$\theta_t:=|o_t|+|x_t|$ and then $\theta_t \to 0$ \almostsurely. Then exists 
$a_t \to 0$, deterministicaly, such that  $\theta_t / a_t \to 0$ \almostsurely. 
From here it follows $|o_t|/a_t\to 0$ and $|x_t|/b_t \to 0$ \almostsurely. 
\end{remark}
We define the stopping times
  \begin{equation}\label{04_1025}
  \tau_R = \inf\{ t : |o_t| \ge R |a_t| \}, \quad \sigma_R = \inf\{ t : |x_t| \ge R |a_t| \}
  \end{equation}
for $R>0$ and
  \begin{equation}\label{04_1633}
  \nu = \min(\tau_R,\sigma_R)\,.
  \end{equation}
%

%
%

From Lemma~\ref{delyon_at} and from (\ref{01_1800})
we conclude that for each $\epsilon>0$ we can choose $R<\infty$ such that
  \begin{equation}
  \P(\nu=\infty) \ge 1-\epsilon.
  \end{equation}
In this way, with a probability so large as we want
we have a deterministic bound common to $|o_t|$ and $|x_t|$.

\vspace{5mm}

Now, consider the similar  process to the algorithm in (\ref{meqal1}) 
but with deterministic step $\gamma_t = 1/(\E_0 t)$ applied
to the function $\varphi(x)=\alpha x$ ($\alpha$ is
the derivative of $\varphi$ in $x^*$),
  \begin{equation}\label{04_1557}
  z_{t} = z_{t-1} - \frac{1}{\E_0 t} (\alpha z_{t-1} + \xi_t), \quad z_0=x_0.
  \end{equation}
Asymptotical properties of this process are known (for example,  
Nevel'son e Has'minskii \cite{b001}). So
  \begin{eqnarray}
  z_t t^{1/2-\epsilon} \to 0,\textrm{\almostsurely}, \textrm{ for each $\epsilon>0$, } \nonumber  \\
  \E |z_t|^2 \le K/t, \quad K>0 \nonumber \\
  \sqrt t z_t \convdist N(0, V). \label{06_1200}
  \end{eqnarray}
where $V$ is the matrix defined in (\ref{m07_1146}).

Based on Lemma~\ref{lema_convdist} in the reference Section,  
Lemma~\ref{lema4_casomultidim} 
will show that, assimptotically, $\sqrt t x_t$ and $\sqrt t z_t$ 
will have the same limiting distribuition, described in (\ref{06_1200}). 
\hfill$\Box$

\begin{lemma}
Consider the following recursive formula, where $b>0$, $a_0$ are real numbers, 
  \begin{equation}\label{erec}
  0 \le a_{t+1} \le (1-\frac{b}{t}) a_t + \OO((t^{-1}),	\quad t=1,2,\ldots\,.
  \end{equation}
Then $a_t \to 0$.
\end{lemma}
\begin{demo}
Consider the recursive sequence, where $\epsilon$ is a positive real number, 
  \[
  0 \le A_{t+1} \le (1-\frac{b}{t}) A_t + \epsilon/t, \quad t=t_0,t_0+1,\ldots\,.
  \]
Then
  \[
  0 \le A_{t+1} \le A_t - \frac{b A_t-\epsilon}{t}, \quad t=t_0,t_0+1,\ldots\,.
  \]
or
  \[
  0 \le bA_{t+1}-\epsilon \le bA_t-\epsilon - b\frac{b A_t-\epsilon}{t}, 
  \quad t=t_0,t_0+1,\ldots\,.
  \]
We write $B_t= bA_t - \epsilon$ and
  \[
  B_{t+1} = B_t (1-b/t)
  \]
so $B_t \to 0$, therefore $A_t \to \epsilon/b$.

Lemma's sequence is 
  \[
  0 \le a_{t+1} \le (1-\frac{b}{t}) a_t + \OO((1)/t, \quad t=1,2,\ldots\,.
  \]
for which we choose $\epsilon>0$ such that 
$o(1)<\epsilon$ if $t\ge t_0$ for some $t_0$.
We define 
  \[
  A_{t+1} = (1-\frac{b}{t}) A_t + \epsilon/t, \quad t=t_0,t_0+1,\ldots
  \]
and $A_{t_0} = a_{t_0}$. Now, we show $0\le a_t \le A_t$ using an induction argument.
Suppose
$A_t - a_t \ge 0$ for $t\ge t_0$. For $t+1$ 
  \[
  A_{t+1} - a_{t+1} =  (1-\frac{b}{t}) (A_t-a_t) + (\epsilon-o(1))/t
  \]
verifying that $A_{t+1} - a_{t+1} \ge 0$ using hypothesis. Then
$ 0 \le a_t \le A_t$.

With $A_t \to \epsilon/b$ and since we can choose a small enough 
$\epsilon$, 
we conclude that $A_t\to0$ 
and therefore $a_t\to0$.
\end{demo}

\begin{lemma}\label{c2lemamatriz}
Let $A$ be a positive definite matrix and symmetrical, $a$, $b$, $c$ and $d$ 
real vectors. Then
  \begin{eqnarray*}
  (a+b+c+d)^T A (a+b+c+d) & \le &
                       a^TAa +  \\
                   & & + 3(b^TAb+c^TAc+d^TAd) + \\ 
                   & & + a^TAb + b^TAa + \\
                   & & + 2a^TA(c+d)\,.
  \end{eqnarray*}
\end{lemma}
\begin{demo}
From
  \begin{eqnarray*}
  (a-b)^T A (a-b)   &=& a^TAa + b^TAb - a^TAb - b^TAa \ge 0 \Leftrightarrow \\
    & \Leftrightarrow & a^TAb + b^TAa \le a^TAa + b^TAb
  \end{eqnarray*}
we have
  \begin{eqnarray*}
  (a+b)^T A (a+b) &  = & a^TAa + b^TAb + a^TAb + b^TAa \\
                  &\le & a^TAa + b^TAb + a^TAa + b^TAb \\ 
                  &=   & 2(a^TAa + b^TAb)\,.
  \end{eqnarray*}
In a similar way 
  \begin{eqnarray*}
  (a+b+c)^T A (a+b+c) & = & a^TAa + b^TAb + c^TAc + \\
                      &   & (a^TAb + b^TAa) + (a^TAc + c^TAa) + \\
                      &   & (b^TAc + c^TAb) \\
      & \le & a^TAa + b^TAb + c^TAc + \\
      &     & (a^TAa + b^TAb) + (a^TAa + c^TAc) + \\
      &     & (b^TAb + c^TAc)  \\
      &   = & 3(a^TAa + b^TAb + c^TAc)\,.\\
  \end{eqnarray*}
So,
  \begin{eqnarray*}
  (a+b+c+d)^T A (a+b+c+d) 
     & =   & (a+(b+c+d))^T A (a+(b+c+d)) \\
     & =   & a^TAa + 
             a^T A (b+c+d) + \\
     &     & (b+c+d)^T A a + 
             (b+c+d)^T A (b+c+d) \\
     & \le & a^TAa + 
             3(b^TAb+c^TAc+d^TAd) +\\
     &     & a^TAb + b^TAa + 
             2 a^T A (c+d) \,.
  \end{eqnarray*}
\end{demo}

\begin{lemma}\label{lema4_casomultidim}
Let $\Delta_t := x_t - z_t$. Then $\sqrt t \Delta_t \convprob 0$.
\end{lemma}
\begin{demo}
From Lemma~\ref{lema_st_draft}, $\gamma_t = \frac{1}{s_t} =\frac{1}{E_0 t}(1+o_t)$ where
$o_t$ is a random variable of ${\cal F}_t$ which converges to $0$ \almostsurely.
Then, from (\ref{meqal1}), (\ref{meqal2}) with $\gamma_t=1/s_t$,  
 \begin{equation}
 x_{t+1} = x_t - \frac{1}{E_0 t}(1+o_t)(\varphi(x_t)+\xi_{t+1})
 \end{equation}
and
\[
x_{t+1}  =  x_t   -   \frac{1}{E_0 t} \varphi(x_{t}) - 
                      \frac{1}{E_0 t} \xi_{t+1} - 
                      \frac{o_t}{E_0 t} \varphi(x_{t}) - 
                      \frac{o_t}{E_0 t} \xi_{t+1}\,.
\]
From Assumption~B3.4, 
  \[
  \varphi(x) = (\varphi(x) - \varphi'(0) x) + \varphi'(0) x\,,
  \]
so
\begin{eqnarray*}
x_{t+1}  =  x_t & - & \frac{1}{E_0 t} \varphi'(0)x_{t} - 
                      \frac{1}{E_0 t} \xi_{t+1} - 
                      \frac{o_t}{E_0 t} \xi_{t+1} - \\
                & - & \frac{1}{E_0 t} 
                        \left(
                          o_t \varphi(x_{t}) + \varphi(x_{t}) - \varphi'(0)x_{t}
                        \right)\,.
\end{eqnarray*}
Define
\[
 v_t  := o_t \frac{\varphi(x_t)}{|x_t|} + 
             \frac{\varphi(x_t) - \varphi'(0)x_t}{|x_t|}
\]
and for $t\le \nu$ we have $|x_t| \le R a_t$ and $|o_t| \le R a_t$ 
\begin{eqnarray}
|v_t| & \le & R a_t \sup_x \frac{|\varphi(x)|}{|x|} + 
                   \sup_{|x|\le Ra_t} \frac{|\varphi(x_t)-\varphi'(0)x_t|}{|x_t|} \le \nonumber \\
      & \le & R a_t M + o(1) := c_t\,. 
\end{eqnarray}
We note that $c_t \to 0$ where $c_t$ is a positive decreasing sequence
and 
\[
x_{t+1}  =  x_t   -   \frac{1}{E_0 t} \varphi'(0)x_{t} -
                      \frac{1}{E_0 t} \xi_{t+1} - 
                      \frac{o_t}{E_0 t} \xi_{t+1} -
                      \frac{1}{E_0 t} v_t |x_t|\,.
\]
Considering the algorithm for $z_t$
  \begin{eqnarray*}
  z_{t+1} & =  z_t & - \frac{1}{E_0 t} (\varphi'(0) z_t + \xi_{t+1}) = \\
          & =  z_t & - \frac{1}{E_0 t} \varphi'(0) z_t - \frac{1}{E_0 t} \xi_{t+1}
  \end{eqnarray*}
and
\begin{eqnarray*}
x_{t+1}  &=&  x_t  - \frac{1}{E_0 t} \varphi'(0)x_{t} 
                   - \frac{1}{E_0 t} \xi_{t+1}
                   - \frac{o_t}{E_0 t} \xi_{t+1}
                   - \frac{1}{E_0 t} v_t |x_t|, \\
z_{t+1}  &=&  z_t  - \frac{1}{E_0 t} \varphi'(0) z_{t} 
                   - \frac{1}{E_0 t} \xi_{t+1}
\end{eqnarray*}
from where
\[
\Delta_{t+1} = \Delta_t - \frac{1}{E_0 t} \varphi'(0) \Delta_t
                        - \frac{1}{E_0 t} v_t |x_t|  
                        - \frac{o_t}{E_0 t} \xi_{t+1}\,.
\]

We wish to show that  
$\sqrt t \Delta_t = \sqrt t (x_t - z_t) \convprob 0$
and for that porpouse we define  
$V_t := \Delta_t^T A \Delta_t$ where $A$ is a definite positive matrix 
to be specified.

First we show that $\E[t V_t \I(t<\nu)] \to 0$ and by Theorem~\ref{uintegra},
p.~\pageref{uintegra},
follows $\sqrt t (x_t - z_t) \convprob 0$.
So,
  \begin{eqnarray*}
  V_{t+1} & = & \Delta_{t+1}^T A \Delta_{t+1} = \\
          & = & (\Delta_t 
                         - \frac{1}{E_0 t} \varphi'(0) \Delta_t
                         - \frac{1}{E_0 t} v_t |x_t|  
                         - \frac{o_t}{E_0 t} \xi_{t+1})^T \cdot \\
          &   & \cdot A \cdot \\
          &   & (\Delta_t 
                         - \frac{1}{E_0 t} \varphi'(0) \Delta_t
                         - \frac{1}{E_0 t} v_t |x_t|  
                         - \frac{o_t}{E_0 t} \xi_{t+1}) 
  \end{eqnarray*}
or, after transposition,  
  \begin{eqnarray*}
  V_{t+1} & = & \Delta_{t+1}^T A \Delta_{t+1} = \\
          & = & (\Delta^T_t 
                         - \frac{1}{E_0 t} \Delta^T_t \varphi'(0)^T
                         - \frac{1}{E_0 t} v^T_t |x_t|  
                         - \frac{o_t}{E_0 t} \xi^T_{t+1}) \cdot \\
          &   &          \cdot A \cdot \\
          &   &          (\Delta_t 
                         - \frac{1}{E_0 t} \varphi'(0) \Delta_t
                         - \frac{1}{E_0 t} v_t |x_t|  
                         - \frac{o_t}{E_0 t} \xi_{t+1})\,.
  \end{eqnarray*}
To estimate $V_{t+1}$ we use Lemma~\ref{c2lemamatriz} to obtain
  \[
  V_{t+1} \le V_t + B_t + C_t + D_t
  \]
with $B_t$, $C_t$ and $D_t$ to be specified and
Using $\I(t+1<\nu) \le \I(t < \nu)$ we estimate $\E[(t+1)V_{t+1} \I(t+1<\nu)]$ by
  \begin{eqnarray*}
  \E[(t+1)V_{t+1} \I(t+1<\nu)] 
     & \le  & \E[(t+1)V_t\I(t<\nu)] \\
     &      & + \E[(t+1)B_t\I(t<\nu)] \\ 
     &      & + \E[(t+1)C_t\I(t<\nu)]  \\
     &      & + \E[(t+1)D_t\I(t<\nu)]\,. 
  \end{eqnarray*}

Considering times when $t \le \nu$ we have 
 $|x_t| \le R a_t$ and 
 $|o_t| \le R a_t$. 
For $B_t$, considering $t < \nu$,
  \begin{eqnarray*}
    B_t & =   & \frac{3}{E^2_0 t^2} 
                   \left( 
                       \Delta^T_t \varphi'(0)^T A \varphi'(0) \Delta_t +
                       |x_t|^2 v^T_t A v_t +
                       o^2_t \xi^T_{t+1} A \xi_{t+1}
                   \right) \\
        & \le & \frac{3}{E^2_0}\frac{1}{t^2} 
                   \left( 
                        K_1 \cdot V_t + 
                        |v_t|^2 \cdot |x_t|^2 \cdot |A|  +
                        o^2_t |A| |\xi_{t+1}|^2 
                   \right) \\
        & \le & \frac{3}{E^2_0}\frac{1}{t^2} 
                   \left( 
                        K_1 \cdot V_t + 
                        c^2_t \cdot R^2 a^2_t \cdot |A| +
                        R^2 a^2_t \cdot |\xi_{t+1}|^2 \cdot |A|
                   \right) \\
        & \le & \frac{3}{E^2_0}\frac{1}{t^2} 
                   \left( 
                        K_1 \cdot V_t + 
                        o(1) +
                        o(1)  \cdot |\xi_{t+1}|^2
                   \right)
  \end{eqnarray*}
where $K_1$ is a positive constant such that
  \[ 
  \Delta^T_t \varphi'(0)^T A \varphi'(0) \Delta_t \le 
  K_1 \Delta^T_t A \Delta_t = 
  K_1 V_t.
  \]
%
%
%
From 
  \[
    (t+1)B_t
         \le  \frac{3(t+1)}{E^2_0}\frac{1}{t^2} 
                   \left( 
                        K_1 \cdot V_t + 
                        o(1) +
                        o(1)  \cdot |\xi_{t+1}|^2
                   \right)
  \]
and using
\begin{itemize}
\item $\frac{3(t+1)}{E^2_0}\frac{1}{t^2} \le \frac{K_3}{t}$, for some positive constant $K_3$;
\item $\frac{3(t+1)}{E^2_0}\frac{1}{t^2} o(1)=o(t^{-1})$;
\item $\E[|\xi_{t+1}|^2 ] = tr(S_\xi)$;
\end{itemize}
we have
  \[
   \E[(t+1)B_t\I(t\le \nu)] = \frac{K_3}{t} V_t + o(t^{-1})\,.
  \]

Now we expand $C_t$,
  \begin{eqnarray*}
    C_t & = & \Delta^T_t A  \frac{-1}{E_0 t} \varphi'(0) \Delta_t + 
              \frac{-1}{E_0 t} \Delta^T_t \varphi'(0) A \Delta_t = \\
        & = & \frac{-1}{t} \Delta^T_t ( A \varphi'(0)/E_0 + \varphi'(0)^T/E_0 A) \Delta_t\,.
  \end{eqnarray*}
Aiming and estimate of $C_t$ in a useful way we find a matrix $A$ which verifies 
$A \varphi'(0)/E_0 + \varphi'(0)^T/E_0 A=I+A$ and we use also $I+A \ge (1+\beta)A$ for  
a real positive constant $\beta$.
We write, for $A=A^T$,
  \begin{eqnarray*}
  A \varphi'(0)/E_0 + \varphi'(0)^T/E_0 A & = & I+A  \Leftrightarrow \\
  \varphi'(0)^T/E_0 A + A \varphi'(0)/E_0 & = & I+A  \Leftrightarrow \\
  \varphi'(0)^T/E_0 A - \frac{A}{2} + A \varphi'(0)/E_0 -\frac{A}{2} & = & I \Leftrightarrow  \\
  (\varphi'(0)^T/E_0 - \frac{I}{2})A + A(\varphi'(0)/E_0-\frac{I}{2}) & = & I
  \end{eqnarray*}
and for use Lyapunov's result (Theorem~\ref{theorem:lyapunov}) we write
the last equality as

 \[
    (\frac{I}{2} - \varphi'(0)^T/E_0) A + A (\frac{I}{2} - \varphi'(0)/E_0)  = - I
 \]
where, from Assumption~B3.3, $\frac{I}{2} - \varphi'(0)/E_0$ is negative definite, therefore
solution $A$ exists and is positive definite. Finalizing,  
  \begin{eqnarray*}
    C_t & = & \frac{-1}{t} \Delta^T_t ( A \varphi'(0)/E_0 + \varphi'(0)^T/E_0 A) \Delta_t \\
        & = & \frac{-1}{t} \Delta^T_t (A + I) \Delta_t \\
        & \le & - (1+\beta)\frac{1}{t} V_t
  \end{eqnarray*}

We estimate the last term $D_t$
  \[
    D_t = \frac{-1}{E_0 t} 
          ( 2 \Delta^T_t A v_t \cdot |x_t| + 2 \Delta^T_t A o_t \xi_{t+1} )\,.
  \]
Recall that we are considering $t<\nu$ and
because we can't use $|\Delta_t| \le V_t$ we follow this
\begin{itemize}
\item $x_t = \Delta_t + z_t$ from where $|x_t|^2 \le |\Delta_t|^2 + |z_t|^2$;
\item $2|\Delta_t|^2 \le K_2 V_t$ ($2$ by convenience)
for a certain positive constant $K_2$.
\end{itemize}
Then,
  \begin{eqnarray*}
    2 \Delta^T_t A v_t \cdot |x_t|  
        & \le & 2 |\Delta_t| \cdot |x_t|  \cdot |A| \cdot c_t \\
        & \le & (|\Delta_t|^2 + |x_t|^2)  \cdot |A| \cdot c_t \\
        & \le & (2|\Delta_t|^2 + |z_t|^2) \cdot |A| \cdot c_t \\
        & \le & (K_2 V_t + |z_t|^2) \cdot |A| \cdot c_t
  \end{eqnarray*}        
We considering again the estimation of $D_t$
  \begin{eqnarray*}
  D_t & \le & \frac{-1}{E_0 t} 
          ( 2 \Delta^T_t A v_t \cdot |x_t| + 2 \Delta^T_t A o_t \xi_{t+1} ) \le \\
      & \le & \frac{K_2}{E_0 t}\cdot|A|\cdot c_t \cdot V_t 
              + \frac{1}{E_0 t}\cdot|A|\cdot c_t \cdot |z_t|^2
              - \frac{2}{E_0 t}  \Delta^T_t A o_t \xi_{t+1} \,. 
  \end{eqnarray*}
Taking 
\begin{itemize}  
\item $\E[ |z_t|^2 ] = K_4/t$, for some constant $K_4$;
\end{itemize}
Then
  \begin{eqnarray*}
  \E[ (t+1) D_t ] & = &
      \frac{K_2(t+1)}{E_0 t}\cdot|A| \cdot c_t \cdot V_t \\
      &     & + \frac{t+1}{E_0 t}\cdot|A|\cdot c_t \cdot \frac{K_4}{t} \\
      & \le & o(1) V_t + o(t^{-1})
  \end{eqnarray*}

Now, putting all together, always considering $t<\nu$,
  \begin{eqnarray*}
  (t+1)V_{t+1} & \le & (t+1) V_t + \frac{K_3}{t} V_t + \\
               &     & o(t^{-1}) - \frac{t+1}{t}(1+\beta) V_t + \\
               &     & o(1) V_t + o(t^{-1}) \le \\
               & \le & V_t (t+1 \frac{K_3}{t} - (1+\beta)\frac{t+1}{t} + o(1)) + o(t^{-1}) \le \\
               & \le & t\cdot V_t (1 + \frac{1}{t} + 
                       \frac{K_3}{t^2} -(1+\beta)\frac{t+1}{t^2} + o(t^{-2})) + o(t^{-1})  \le \\
               & \le & t V_t (1 -(1+\beta)\frac{1}{t} + o(t^{-1})) + o(t^{-1})  \le \\
               & \le & t V_t (1 - (1+\beta+o(1))\frac{1}{t} ) + o(t^{-1}) \le \\
               & \le & t V_t (1 - (\beta/2)\frac{1}{t} ) + o(t^{-1})\,.
  \end{eqnarray*}
It follows that, 
  \[
  \E[(t+1) V_{t+1} \I(t+1 < \nu)] \le \E[t V_{t} \I(t < \nu)] + o(t^{-1})
  \]
and by Lemma~\ref{c2lemamatriz}
  \[
  \E[t V_{t} \I(t < \nu)]  \to 0,
  \]
then, by Theorem~\ref{uintegra}, 
  \[
  t V_{t} \I(t < \nu)  \convprob 0,
  \]
or
  \[
  \sqrt t (x_t - z_t) \I(t < \nu) \convprob 0\,,
  \]
or even, by definition of convergence in probability,
  \[
   \forall \eta>0\quad \P( |\sqrt t (x_t - z_t) \I(t < \nu)| < \eta) \to 1\,.
  \]
The following events are related by 
  \[
  \sqrt t (x_t - z_t) < \eta \Rightarrow \sqrt t (x_t - z_t) \I(t < \nu) < \eta
  \]
and by $P(\sqrt t (x_t - z_t) < \eta) \le \P(\sqrt t (x_t - z_t) \I(t < \nu) < \eta)$ we have
  \[
  \sqrt t (x_t - z_t) \convprob 0\,.
  \]
 
\end{demo}

\section{Some standard results}\label{resultadosusados}

\begin{theorem}[A. M. Lyapunov, 1947 (cited in \cite{b002}, Chap.~13.1)]
\label{theorem:lyapunov}
Let $U,W\in\mathbb{C}^{n\times n}$ and let $W$ be positive definite.
\begin{enumerate}
\item[(a)] If $U$ is stable then the equation
  \[
  UA + AU^* = W
  \]
as a unique solution $A$ negavtive definite.
\item[(b)] If exists a negative definite matrix $A$ satisfying
the above equation then $A$ is stable.
\end{enumerate}
\end{theorem}
\begin{remark}
Stablçe is when all eigenvalues are negative.
When all eigenvalues are negative  then the matrix is negative definite.
\end{remark}


\begin{lemma}[Markov Inequality (for example, \cite{b016})]\label{markovineq}
Let $Z$ a r.v. and $g:\mathbb{R} \to [0,\infty]$ a non decreasing function. 
Then
  \[
   \E g(Z) \ge \E( g(Z); Z \ge c ) \ge g(c) \P(Z \ge c)
  \]
\end{lemma}


\begin{theorem}[Martingale convergence, \cite{b016}, Cap. 12]\label{martingale}
Let $M$ be a martingale for which $M_n \in \LL^2, \forall n$.
Then $M$ is limitied in $\LL^2$ iif
\[
  \sum \E[ (M_k - M_{k-1})^2 ]  < \infty
\]
and when this we have
\[
  M_n \to M_\infty \textrm{ \almostsurely\ and in } \LL^2\,.
\]
\end{theorem}


\begin{theorem}[\cite{b016}, Chap. 13.7]\label{uintegra}
Let $(X_n)$ be a sequence in $\LL^1$ and $X \in \LL^1$. 
Then $X_n \to X$ in $\LL^1$, or similarly
$\E(|X_n - X|) \to 0$, iif, the following conditions are verifyed,
\begin{enumerate}
\item $X_n \to X$ in probability;
\item the sequence $(X_n)$ is uniformly integrable ($\forall \epsilon>0 \exists K\,:\,\E[|X|;|X|>K]<\epsilon$).
\end{enumerate}
\end{theorem}


\begin{lemma}[Slutsky's Theorem, \cite{b029} Sec.8.6]\label{lema_convdist}
If $|X_t-Z_t|\convprob 0$ and $X_t$ converges in distribution 
then $Z_t$ converges in distribuition for the same limit.
\end{lemma}


\begin{theorem}[Kolmogorov Law of Iterated Logarithm \cite{b016}] \label{iteratedlaw}
Let $X_1, X_2, \ldots$ be random variables independent and identically 
distributed with mean $0$ and variance 1. Let $S_n := X_1 + \cdots + X_n$. Then,
\almostsurely,
  \[
  \limsup \frac{S_n}{\sqrt{2n\log\log n}}\to+1, \quad\quad \liminf \frac{S_n}{\sqrt{2n\log\log n}}\to -1\,.
  \]
\end{theorem}


\addcontentsline{toc}{section}{\numberline{}Bibliography}
\bibliography{sabib}
\bibliographystyle{plain}

\end{document}